\documentclass[11pt]{article}
\usepackage{amssymb}
\usepackage{color, amsmath,amssymb, amsfonts}
\usepackage{amstext,amsthm, latexsym}

\usepackage{amssymb, epsfig, amssymb, latexsym}
\usepackage{amsmath}
\usepackage{graphicx}
\usepackage{longtable}
\usepackage{subfigure}

  \newtheorem{Main Results}[theorem]{MainResults}

\newcommand{\x}{\mathbf x}
\newcommand{\z}{\mathbf z}
\newcommand{\w}{\mathbf w}
\newcommand{\vv}{\mathbf v}
\newcommand{\F}{\mathbf F}
\newcommand{\HH}{\mathbf H}
\newcommand{\p}{\mathbf P}

\newcommand{\K}{\mathbf K}
\newcommand{\R}{\mathbf R}

\newcommand{\Q}{\mathbf Q}
\begin{document}

\title{State estimation under non-Gaussian L\'evy noise:\\ A modified Kalman filtering method
\footnote{This work was partly supported by the NSF grant DMS-1025422. } }


\author{Xu Sun$^{1,2}$, Jinqiao Duan$^{2,3}$, Xiaofan Li$^{2}$  and Xiangjun Wang$^{1}$
\\
1. Huazhong University of Science and Technology\\
 Wuhan 430074, Hubei, China\\
    E-mail: xjwang@hust.edu.cn\\
2. Illinois Institute of Technology \\
   Chicago, IL 60616, USA \\
   E-mail: lix@iit.edu\\
3. Institute for Pure and Applied Mathematics\\
University of California\\
Los Angeles, CA 90095, USA \\
   E-mail: jduan@ipam.ucla.edu
}

\date{\today}
\maketitle

\begin{abstract}
The Kalman filter is extensively used for state estimation for linear systems under Gaussian noise. When non-Gaussian L\'evy noise is present, the conventional Kalman filter may fail to be effective due to the fact that the non-Gaussian L\'evy noise may have infinite variance. A modified Kalman filter for linear systems with non-Gaussian L\'evy noise is devised. It works effectively  with reasonable computational cost. Simulation results are presented to illustrate this non-Gaussian filtering method.

\medskip



\end{abstract}

\maketitle

\pagestyle{plain}

\section{Introduction and   statement of the problem}
\label{intro}

%
%

The Kalman filter, or the Kalman filtering method, provides an efficient way to estimate the state   of a linear dynamical system subject to Gaussian white noise \cite{BrownHwang1997, GrewalAndrews2008, Jaz}. It has been widely used in applications such as target tracking, parameter estimation, control theory, signal processing, and other data assimilation tasks.

The Kalman filter   requires the noise be either Gaussian or with finite variance \cite{BrownHwang1997, GrewalAndrews2008}, and thus it is not applicable  to  linear systems   with non-Gaussian noise of infinite variance.
As non-Gaussian L\'evy noise with infinite variance  exists ubiquitously  \cite{Woyczynski2001, SornetteIde2001}, it is desirable  to study the Kalman filtering problems under L\'evy noise. Very little work has been done for this issue.     Breton and Musiela \cite{Breton1993}  presented a scheme for Kalman filtering with noise of infinite variance, while assuming the contribution of the jumps are exactly known. The filter in \cite{Breton1993} is nonlinear and recursive, and thus may greatly limit its application in practice.  Ahn and Feldman \cite{AhnFeldman1999} proposed to minimize the difference between the true state   and the filtered observation in the $L^\mu$-norm. However, as pointed in \cite{SornetteIde2001}, this method does not really address the Kalman filtering problem which consists of combining forecasts and observations. The method in \cite{SornetteIde2001} focuses on large errors and has a robust performance but  high computational cost due to the matrix diagonalization and the operation of the fractional power in each step. This may also be  an obstacle for real-time implementation in practical applications. Note that in practice, each iteration step must be completed during every sampling period, and it is greatly desirable to make the algorithm as fast as possible.

In this paper, we will  present an algorithm, which has the similar computational cost as that of the Kalman filter, but can be applied to linear systems with non-Gaussian L\'evy noise of infinite variance.

We consider the following discrete time model
with the state equation
\begin{align}\label{eq_sec1_1}
\x_{k+1} = \F_k \x_k + \w_k,
\end{align}
and the observation equation
\begin{align}\label{eq_sec1_2}
\z_k=\HH_k \x _k + \vv_k,
\end{align}
where $\x_k$, an $n$-by-$1$ vector, is the state variable, and $\z_k$, an $m$-by-$1$ vector, is the measurement (or observation) variable,   $\w_k$ represents the modeling error noise, $\vv_k$ the measurement error noise, and $\F_k$ and $\HH_k$ are $n$-by-$n$ and $m$-by-$n$ matrices, respectively. We only consider the cases where $\w_k$ is a Gaussian noise, and $\vv_k$ is a non-Gaussian L\'evy noise.


This paper is arranged as follows. In section 2, the usual Kalman filter is briefly reviewed. The proposed modified Kalman filter is presented in section 3. A simulation example is provided in Section 4  to illustrate the effectiveness of the modified Kalman filter.

\section{Review of the conventional Kalman  filter}

A derivation of the Kalman filter is briefly reviewed in this section. Some equations and ideas presented in this section will be used to present our   proposed modified Kalman filter  in the next section. Derivations of the Kalman filter can be found in many references   \cite{BrownHwang1997, GrewalAndrews2008, Jaz}.


Consider the model as given in (\ref{eq_sec1_1}) and (\ref{eq_sec1_2}). The Kalman filtering assumes that both the modeling error noise $\w_k$ and the measurement disturbance $\vv_k$ are Gaussian with the following covariance matrix,
\begin{align}\label{eq_sec2_1}
E\left[\w_i \w_k^T\right]=\begin{cases} \mathbf Q_k, &\text{for} \quad i=k,\\0, &\text{for} \quad i\ne k.\\ \end{cases}
\end{align}
\begin{align}\label{eq_sec2_2}
E\left[\vv_i \vv_k^T\right]=\begin{cases} \mathbf R_k, &\text{for} \quad i=k,\\0, &\text{for} \quad i\ne k.\\ \end{cases}
\end{align}

Let $\bar \x_k$ be the priori estimate, which is the estimate of $\x_k$ given $\z_0$, $\z_1$, $\cdots$, $\z_{k-1}$, and let $\hat \x_k$ be the posterior estimate, which is the estimate of $\x_k$ given $\z_0$, $\z_1$, $\cdots$, $\z_{k}$. It is known that
\begin{align}\label{eq_sec2_3a}
E\{\bar \x_k\} = E\{ \x_k\}
\end{align}
and
\begin{align}\label{eq_sec2_3}
\bar \x_{k+1} =\F_k \hat \x_k,
\end{align}
where $E\{\cdot\}$ represents expectation and $\F_k$ is from (\ref{eq_sec1_1}).

\begin{figure}
  \epsfig{file=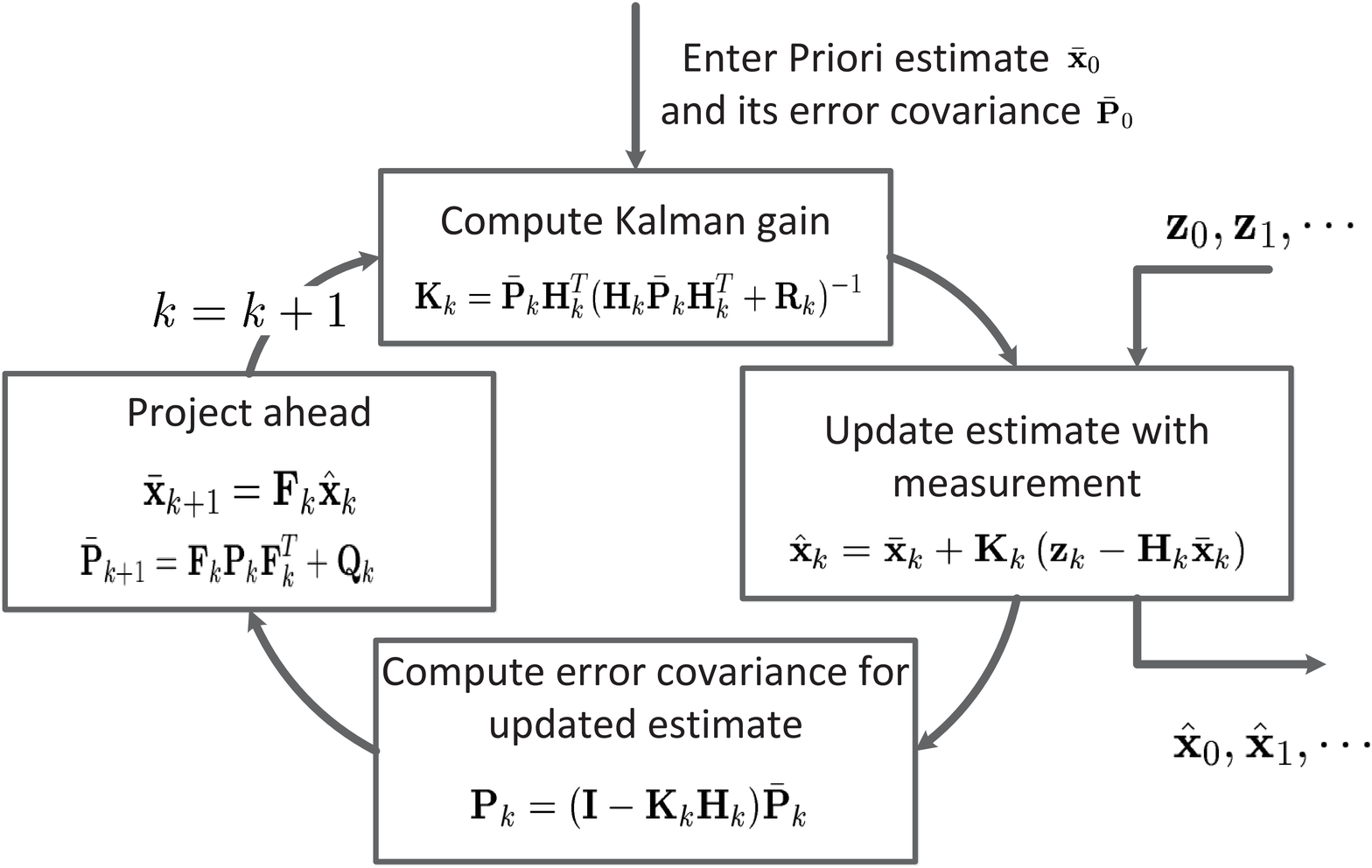,width=\linewidth}
  \caption{The usual Kalman filtering algorithm}
  \label{fig1}
\end{figure}

The Kalman filter assumes that the posterior estimate is expressed as the prior estimate corrected by the measurement data,
\begin{align}\label{eq_sec2_4}
\hat \x_k = \bar \x _k  + \mathbf K _k \left( \z_k-\HH_k \bar \x_k\right),
\end{align}
for some $n$-by-$m$ matrix $\K_k$ (so called Kalman gain). Note that   $\HH_k$ is from (\ref{eq_sec1_2}).
The Kalman gain $\K_k$ is solved by minimizing       $E\left[(\hat \x_k - \x_k)^2\right]$.
Note that
\begin{align}\label{eq_sec2_5}
E\left[(\hat \x_k - \x_k)^T (\hat \x_k - \x_k)\right] = \text{Tr} \{ \p_k\},
\end{align}
where $\text{Tr}\{\cdot\}$ represents the trace operator, and the $n$-by-$n$ covariance matrix $\p_k$ is defined as follows
\begin{align}\label{eq_sec2_6}
 \p_k = E\left[(\x _k -\hat \x_k)(\x _k -\hat \x_k)^T \right].
\end{align}
Define
\begin{align}\label{eq_sec2_7}
\bar \p_k = E\left[(\x _k -\bar \x_k)(\x _k -\bar \x_k)^T \right].
\end{align}
It follows from (\ref{eq_sec2_1}), (\ref{eq_sec2_2}), (\ref{eq_sec2_6})  and (\ref{eq_sec2_7}) that
\begin{align}\label{eq_sec2_8}
\bar \p_{k+1} = \F_k \p_k \F_k ^T + \mathbf Q_k.
\end{align}
Substituting (\ref{eq_sec2_4}) into (\ref{eq_sec2_6}), we   get
\begin{align}\label{eq_sec2_9}
\p_k = (\mathbf I -\mathbf K_k \HH_k) \bar \p_k (\mathbf I -\mathbf K_k \HH_k)^T + \mathbf K_k \mathbf R_k \mathbf K_k^T.
\end{align}
It follows from (\ref{eq_sec2_9}) that
\begin{align}\label{eq_sec2_10}
\frac{d}{d\mathbf K_k}Tr\{\p_k\} = -2(\HH-K \bar\p_k)^T +2 \K_k(\HH_k \bar \p_k \HH_k^T + \R_k).
\end{align}
Solve $\K_k$ by letting $\frac{d}{d\mathbf K_k}Tr\{\p_k\} = 0$, we get
\begin{align}\label{eq_sec2_11}
\K_k = \bar \p_k \HH_k^T(\HH_k \bar \p_k \HH_k^T + \R_k)^{-1}
\end{align}
By   (\ref{eq_sec2_9}) and (\ref{eq_sec2_11}), $\p_k$ can be rewritten as
\begin{align}\label{eq_sec2_12}
\p_k = (\mathbf I- \K_k \HH_k) \bar \p_k.
\end{align}
Combining (\ref{eq_sec2_3}), (\ref{eq_sec2_11})  and (\ref{eq_sec2_12}), we thus have the conventional Kalman filter. This algorithm   is shown in Figure \ref{fig1}.

\section{A modified Kalman filter}

It is   known that the discrete time Gaussian white noise can be approximated by the increments of  Brownian motion  per time step, and the non-Gaussian L\'evy noise can be approximated by the increments of the corresponding L\'evy process per time step. By L\'evy-Ito theorem \cite{Applebaum2009}, a L\'evy process can be decomposed into the sum of a Gaussian process and a pure jump process. It is shown in \cite{AsmussenRosinski2001} that the small jumps of a L\'evy process can be approximated by a Gaussian process. Therefore, we can approximately regard a L\'evy process as combination of a Gaussian process and a process with big jumps. For more information about decomposition of a L\'evy processes, see \cite{AsmussenRosinski2001, Applebaum2009}. These results enable us to decompose a non-Gaussian L\'evy noise into a Gaussian white noise plus some extremely large values.

In our proposed filtering method,  we convert the original L\'evy noise into a Gaussian white noise by clipping off its extremely large values.

Let $\tilde \vv_k$ represent the clipped version of the L\'evy measurement disturbance $\vv_k$, and let $\tilde \z_k$ represent the corresponding clipped observation. Thus
\begin{align}\label{eq_sec3_1}
\tilde \z_k = \HH_k \x_k + \tilde \vv_k.
\end{align}
In practice, since the measurement noise, $\vv_k$, is unknown, we propose to clip the observation  $\z_k$ instead of $\vv_k$ in a component-wise way by the following operation:
\begin{align}\label{eq_sec3_2}
\tilde \z_k^i=  \begin{cases} \sum_j \HH_k^{i,j} \bar \x_k^j + C \cdot \text{sign} \left(\z_k^i-\sum_j \HH_k^{i,j} \bar \x_k^j \right)   & \text{if} \quad |\z_k^i -\sum_j \HH_k^{i,j} \bar \x_k^j |\ge C,\\
 \z_k^i  & \text{if} \quad |\z_k^i -\sum_j \HH_k^{i,j} \bar \x_k^j | < C,\\
\end{cases}
\end{align}
where $C$ is some positive threshold value, $\z_k^i$ and $\bar \x_k^i$ represent the $i$-th components of the vectors $\z_k$ and $\bar \x_k$, respectively, and $\sum_j \HH_k^{i,j} \bar \x_k^j$ is the $i$-th component of the vector $\HH_k \bar \x_k$. Note that $C$  is determined by the statistical properties of the measurement noise $\vv_k$.
Replacing the observation value $\z_k$ in (\ref{eq_sec2_4}) with its clipped value, we get
\begin{align}\label{eq_sec3_3}
\hat \x_k = \bar \x _k  + \mathbf K _k \left( \tilde \z_k-\HH_k \bar \x_k\right).
\end{align}
Repeating the same procedure in Section  $2$ to solve the Kalman gain $\K_k$ by minimizing $E\{(\x_k -\hat \x_k)^2\}$, we get
\begin{align}\label{eq_sec3_4}
\K_k = \bar \p_k \HH_k^T(\HH_k \bar \p_k \HH_k^T + \tilde \R_k)^{-1},
\end{align}
where
$\tilde \R_k$ is the covariance matrix of $\tilde \vv_k$ defined as
\begin{align}\label{eq_sec3_5}
\tilde \R_k = E\left \{ \tilde \vv_k  \tilde \vv_k  ^T\right\}.
\end{align}
In the conventional Kalman filter, $\Q_k$ and $\R_k$ are assumed to be known, and as noted in \cite{GrewalAndrews2008}, it is often a difficult task to estimate the covariance matrices $\Q_k$ and $\R_k$.

In the modified Kalman filter here, we only assume $\Q_k$ is known and suggest $\tilde \R_k$ be estimated as follows.
It follows from (\ref{eq_sec3_1}) and (\ref{eq_sec3_2}) that
\begin{align}\label{eq_sec3_6}
\tilde \R_k =&E\left \{ \tilde \vv_k  \tilde \vv_k  ^T\right\}\nonumber\\
 =&E\left \{\left[(\tilde \z_k - \HH_k \bar \x_k) - \HH_k (\x_k -\bar \x_k) \right]\left[(\tilde \z_k - \HH_k \bar \x_k) - \HH_k (\x_k -\bar \x_k) \right]^T\right\}\nonumber\\
 =& (\tilde \z_k - \HH_k \bar \x_k) (\tilde \z_k - \HH_k \bar \x_k) ^T + \HH_k \bar \p_k \HH_k^T.
\end{align}
In deriving the last identity of (\ref{eq_sec3_6}), we have used the fact that $\tilde \z_k$ and $\bar \x_k$  are known values and
\begin{align}\label{eq_sec3_7}
E\left \{ (\tilde \z_k - \HH_k \bar \x_k)(\x_k -\bar \x_k)^T \HH_k ^T \right\}
 = (\tilde \z_k - \HH_k \bar \x_k) E\{ (\x_k -\bar \x_k)  ^T\}\HH_k^T=0.
\end{align}
With (\ref{eq_sec3_6}), (\ref{eq_sec3_4}) can be rewritten as
\begin{align}\label{eq_sec3_8}
\K_k = \bar \p_k \HH_k^T(2\cdot \HH_k \bar \p_k \HH_k^T + \tilde {\tilde {\R}}_k)^{-1},
\end{align}
where
\begin{align}\label{eq_sec3_9}
\tilde{\tilde {\R}}_k = (\tilde \z_k - \HH_k \bar \x_k) (\tilde \z_k - \HH_k \bar \x_k) ^T.
\end{align}

\begin{figure}
  \epsfig{file=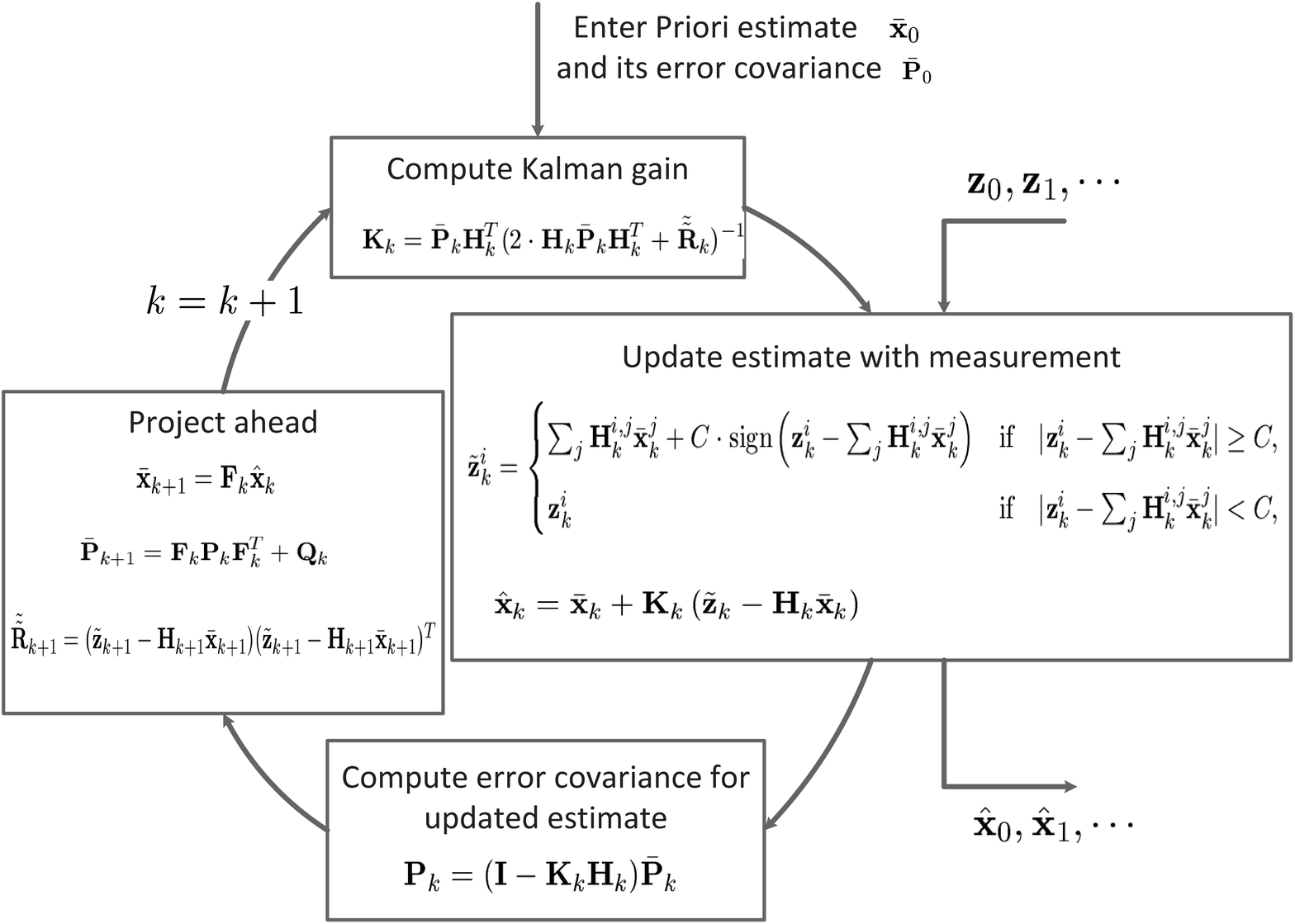,width=\linewidth}
  \caption{The modified Kalman filtering algorithm}
  \label{fig2}
\end{figure}

Combining   equations (\ref{eq_sec3_2}), (\ref{eq_sec3_3}), (\ref{eq_sec3_8}), and (\ref{eq_sec3_9}), we   obtain the modified Kalman filter, as shown graphically in Figure \ref{fig2}.

Comparing with the conventional Kalman filter, the proposed filter has an moderately increased computational cost due to the following two operations: i) the clipping operation for $\z_k$; ii) the computation of $\tilde {\tilde {R}}_k$. The former operation is implemented by IF-ElSE sentence, and the latter is simply a vector-vector outer product.


\section{Simulation results}

 Consider a particle moving in the plane at some velocity subject to random perturbations in its trajectory. The new position at time $k+1$ is equal to the old position at time $k$ plus the velocity and noise. The model can be expressed in form of     (\ref{eq_sec1_1}) and (\ref{eq_sec1_2}) as
\begin{align}\label{example_1}
\begin{pmatrix} x^1_{k+1}\\ x^2_{k+1}\\ u^1_{k+1}\\ u^2_{k+1} \end{pmatrix}= \begin{pmatrix} &1 & 0&1& 0\\ &0 &1 &0 &1\\ &0 &0 &1 &0\\ &0 &0 &0 &1\end{pmatrix} \begin{pmatrix} x^1_{k}\\ x^2_{k}\\ u^1_{k}\\ u^2_{k} \end{pmatrix} + \begin{pmatrix}   w^1_k\\  w^2_k\\  w^3_k \\w^4_k \end{pmatrix},
\end{align}
and
\begin{align}\label{example_2}
\z_k = \begin{pmatrix} &1 &0 & 0& 0 \\ &0 &1 & 0& 0 \end{pmatrix} \x_k + \begin{pmatrix} v^1_k\\ v^2_k\end{pmatrix},
\end{align}
where ($x^1_k$, $x^2_k$) is the position at time $k$,  $u^1_k$, $u^2_k$  the velocity, $w^1_k$, $w_k^2$, $w_k^3$, and $w_k^4$ are all Gaussian white noises with zero mean and unit variance, and $v_k^1$ and $v_k^2$ are independent and identically distributed  noises consist of two components: i) a symmetric $\alpha$-stable L\'evy noises with the index of stability $\alpha=1.3$ and the scale parameter $\sigma = 10$ (see \cite{JW}); ii) a Gaussian white noise with variance of 5.  Since the measurement noises, $v_k^1$ and $v_k^2$, have   infinite variances, the conventional Kalman filter can not be applied to estimate the position  $x_k^1$ and $x_k^2$. So we apply the modified Kalman filter proposed in the previous section.

\begin{figure}
  \epsfig{file=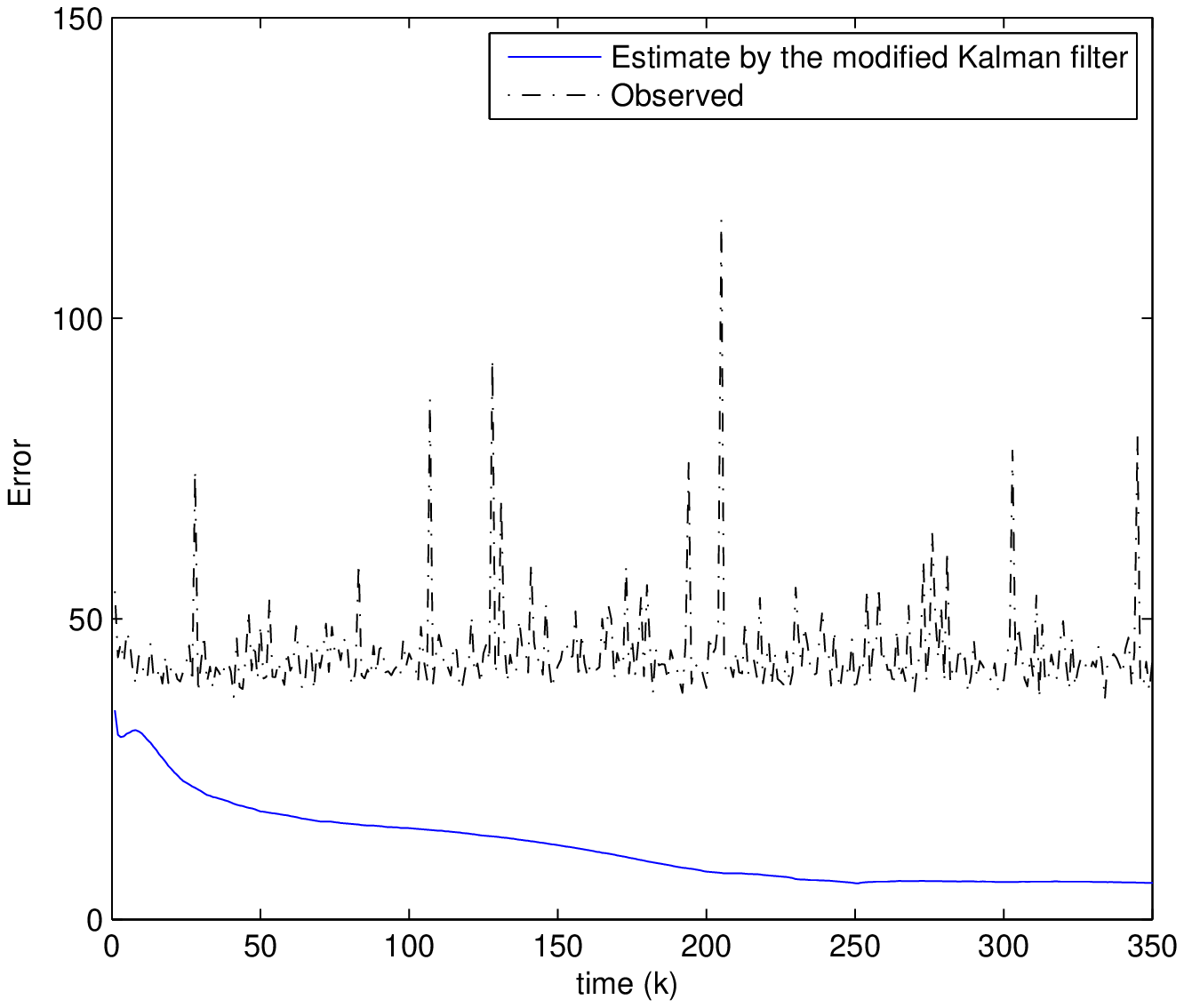,  width=\linewidth}
  \caption{The error of the modified Kalman filtering algorithm}
  \label{fig3}
\end{figure}

Take $x^1_0=10$, $x^2_0=10$, $u^1_0=1$,  $u^2_0 = 0$, and we apply the   modified Kalman filtering method to estimate $x_k^1$ and $x_k^2$. In the simulation, the initial a priori estimate of the state, ($\bar x^1_0$, $\bar  x^2_0$), is set to be equal to the observation at time $0$, its error covariance, $\bar P_0$, is set to be unit matrix, and the threshold value $C$ is set to be $40$. The simulation results  are shown in Figure \ref{fig3}, where the estimate position error, $ER$,   defined by
\begin{align}
ER_k =\sqrt{( z^1_k -  x^1_k)^2 +( z^2_k -  x^2_k)^2 },
\end{align}
is compared with the observed position error, $OR$, defined by
\begin{align}
OR_k =\sqrt{(\bar x^1_k -  x^1_k)^2 +( \bar x^2_k - x^2_k)^2 }.
\end{align}
 The results in Figure \ref{fig3} are calculated by averaging $10,000$ times of simulations.
It is seen from this figure that the position estimation error is significantly improved by using our modified Kalman filter.

In the simulation, we select $C$ by the method of trial and error, and it is found that the threshold value $C$ is not very picky and $C$ can vary from $30$ to $100$ without significant effects on the performance of the modified Kalman filter. Determining the optimal $C$, which is crucial for the modified Kalman filter, deserves further research and will be left for our future work.


%

\newpage

\newpage

\newpage


\end{document}